\documentclass[12pt]{amsart}
\usepackage{amsfonts,amssymb,graphicx, enumerate,youngtab}

\input prepictex
\input pictex
\input postpictex

\numberwithin{equation}{section}

\newtheorem{theorem}{Theorem}[section]

\theoremstyle{definition}

\theoremstyle{remark}

\newcommand{\CC}{\mathbb{C}}

\newcommand{\ZZ}{\mathbb{Z}}

\newcommand{\hh}{\mathfrak{h}}

\newcommand{\la}{\langle}
\newcommand{\ra}{\rangle}

\begin{document}
\title{Macdonald polynomials as characters of Cherednik algebra modules}

\author{Stephen Griffeth}
\address{Stephen Griffeth \\
Instituto de Matem\'atica y F\'isica \\
Universidad de Talca \\
sgriffeth@inst-mat.utalca.cl}
\email{}

\thanks{I am extremely grateful to Adriano Garsia and Luc Lapointe for very enlightening discussions, which were indispensable for the discovery of these results. I am also grateful to Iain Gordon for helpful comments on an early version of this paper, and to Luc Lapointe for empirical verification of the main result. During the time this research was carried out I have benefitted from the financial support of Fondecyt Proyecto Regular 1110072.}

\begin{abstract}
We prove that Macdonald polynomials are characters of irreducible Cherednik algebra modules. 
\end{abstract}

\maketitle

\section{Introduction}

The \emph{Macdonald positivity theorem} \cite{Hai} asserts that the coefficients $\tilde{K}_{\lambda,\mu}(q,t)$ defined by the equation
$$\tilde{H}_\mu=\sum \tilde{K}_{\lambda,\mu}(q,t) s_\lambda$$ are polynomials in $q,t$ with non-negative coefficients. Here $\tilde{H}_\mu$ is Haiman's \cite{Hai3} renormalized, plethystically transformed Macdonald polynomial $\tilde{H}_\mu$, and $s_\lambda$ is the Schur function.  Several proofs (\cite{As}, \cite{BeFi}, \cite{Gor3}, \cite{GrHa}) have now appeared, using an astonishing range of tools. Haiman \cite{Hai} proved the Macdonald positivity conjecture by identifying $\tilde{H}_\mu$ as the bigraded character of a certain fiber of a rank $n!$ vector bundle---the \emph{Procesi bundle}---on the Hilbert scheme $\mathrm{Hilb}_n(\CC^2)$ of $n$ points in $\CC^2$. The fact that the Procesi bundle is a vector bundle is often called the \emph{$n!$ theorem}. As of this writing, Haiman's original proof is the only one that establishes the $n!$ theorem. Our main result, Theorem~\ref{main}, is that the $n!$ theorem is equivalent to a character formula for irreducible Cherednik algebra modules.

Kashiwara and Rouquier \cite{KaRo} quantized the Hilbert scheme, constructing a sheaf of algebras on it whose algebra of global sections is the spherical rational Cherednik algebra. Iain Gordon's work \cite{Gor2} on diagonal coinvariant rings for real reflection groups specializes in type $A$ to identify a certain quotient of the quantization of the Procesi bundle with the diagonal coinvariant ring (this identification depends on Haiman's formula for the dimension of the diagonal coinvariant ring). Combining our theorem with Haiman's work \cite{Hai} identifies other quotients with Garsia-Haiman modules, which are the fibers of the Procesi bundle at torus-fixed points. 

\section{Garsia-Haiman modules and irreducible representations of the restricted Cherednik algebra}

\subsection{} Let $\lambda$ be a partition of $n$ and write $S^\lambda$ for the corrresponding irreducible $S_n$-module. Let $n(\lambda)$ denote the degree of the unique minimal occurence of $S^\lambda$ in $\CC[x_1,x_2,\dots,x_n]$. Let $I_\lambda \subseteq \CC[x_1,\dots,x_n]$ be the largest homogeneous $S_n$-invariant ideal that does not contain this copy of $S^\lambda$. According to Garsia and Procesi \cite{GaPr}, the graded $S_n$-character of the ring $R_\lambda(x)=\CC[x_1,x_2,\dots,x_n]/I_\lambda$ is the Hall-Littlewood polynomial $t^{n(\lambda)} H_\lambda(z,t^{-1})$. In particular, the only irreducible $S_n$-modules that occur in $R_\lambda(x)$ are of the form $S^\mu$ for $\mu \geq \lambda$ in dominance order.

\subsection{} Let $R_{\lambda}(y)=\CC[y_1,y_2,\dots,y_n]/I_{\lambda}$ and define the \emph{Garsia-Haiman} module $R_\lambda$ (introduced in \cite{GaHa}) to be the unique bigraded Gorenstein quotient of $R_{\lambda^t}(x) \otimes_\CC R_\lambda(y)$ in which the sign representation---which occurs exactly once in  $R_{\lambda^t}(x) \otimes R_\lambda(y)$, in bidegree $(n(\lambda^t),n(\lambda))$---survives (see \cite{Hai2} Proposition 4.1.2 and the paragraph preceding it). Haiman's $n!$ theorem is 
\begin{theorem}[Haiman  \cite{Hai}]
$\mathrm{dim}(R_\lambda)=n!$. 
\end{theorem} Haiman showed in \cite{Hai3} that this theorem implies $\mathrm{ch}(R_\lambda)=\tilde{H}_\lambda$.

\subsection{} Let $H$ be the rational Cherednik algebra for the symmetric group ``at $t=0$". That is, $H$ is the $\CC$-algebra generated by the symmetric group $S_n$, pairwise commuting variables $x_1,x_2, \dots,x_n$, and pairwise commuting variables $y_1,y_2,\dots,y_n$ subject to the additional relations
$$w x_i w^{-1}=x_{w(i)} \quad \text{and} \quad w y_i w^{-1}=y_{w(i)} \quad \hbox{for $1 \leq i \leq n$ and $w \in S_n$,}$$
$$y_i x_i=x_i y_i-\sum_{j \neq i} s_{ij} \quad \hbox{for all $1 \leq i \leq n$,}$$ and
$$y_i x_j=x_j y_i+s_{ij} \quad \hbox{for all $1 \leq i \neq j \leq n$.}$$ In these formulas we have written $s_{ij}$ for the transposition interchanging $i$ and $j$.

\subsection{} The defining relations show that $H$ is graded by assigning degree $0$ to $S_n$, degree $1$ to each $x_i$, and degree $-1$ to each $y_i$. We will write $H^d$ for the degree $d$ piece of this grading. There is also an increasing filtration defined by total degree, 
$$W_l H=\CC \{f(x) g(y) w \ | \ f \in \CC[\hh], \ g \in \CC[\hh^*], \ w \in S_n, \ \mathrm{deg}(f)+\mathrm{deg}(g) \leq l \}.$$ The PBW theorem for $H$ asserts that the associated graded algebra is
$$\mathrm{gr} H= \CC[x_1,\dots,x_n,y_1,\dots,y_n] \rtimes S_n,$$ and since $W_l H$ is graded, $\mathrm{gr} H$ is bigraded with $x$ in degree $(1,1)$ and $y$ in degree $(1,-1)$. By Proposition~3.6 of \cite{Gor}, symmetric polynomials in the $x$ variables are central in $H$, and so are symmetric polynomials in the $y$ variables. Let $\overline{H}$ be the quotient of $H$ by the ideal generated by all positive degree elements of $\CC[\hh]^{S_n}$ and all positive degree elements of $\CC[\hh^*]^{S_n}$. 

\subsection{} According to Proposition 4.3 of \cite{Gor} (which is a an application of general results of Holmes and Nakano \cite{HoNa} to our situation), the irreducible modules of $\overline{H}$ are indexed by partitions $\lambda$ of $n$, and we write $L(\lambda)$ for the unique irreducible generated by a subspace of vectors annihilated by $\hh$ and carrying the representation $S^\lambda$. Corollary 1.14 of \cite{EtGi} implies that, as an $S_n$-module, $L(\lambda) \cong \CC S_n$. Write $t_\lambda \in L(\lambda)$ for a non-zero vector carrying the trivial representation, and define a filtration by $W_l L(\lambda)=W_l H t_\lambda$ for $l \in \ZZ_{\geq 0}$. We shift the grading on $L(\lambda)$ so that $t_\lambda$ lies in degree zero. The spaces $W_l L(\lambda)$ are graded subspaces and hence $\mathrm{gr}L(\lambda)$ is bigraded. Our theorem is
\begin{theorem} \label{main}
The following are equivalent:
\begin{itemize}
\item[(a)] The $n!$ theorem holds;
\item[(b)] the $\CC[x_1,\dots,x_n,y_1,\dots,y_n] \rtimes S_n$-module $\mathrm{gr} L(\lambda)$ is isomorphic to $R_\lambda$;
\item[(c)] the character of $\mathrm{gr} L(\lambda)$ is the Macdonald polynomial $\tilde{H}_\lambda(qt,q t^{-1})$;
\item[(d)] we have $\mathrm{dim}_\CC(\mathrm{gr}^i L(\lambda))=\mathrm{dim}_\CC( \mathrm{gr}^{n(\lambda)+n(\lambda^t)-i} L(\lambda))$ for all $0 \leq i \leq n(\lambda)+n(\lambda^t)$.
\end{itemize}
\end{theorem} The remainder of the paper is devoted to the proof of this theorem.

\subsection{} We first observe that (c) implies (d) follows from the symmetry of Macdonald polynomials given in Corollary 2.7 of \cite{Hai3}. Also, assuming that (b) holds, the dimension of $R_\lambda$ is $n!$, and therefore the results of \cite{Hai3} imply that the bigraded character of $R_\lambda$ is the Macdonald polynomial $\tilde{H}_\lambda(q,t)$. Taking into account the grading shift $\mathrm{deg}(x_i)=(1,1)$ and $\mathrm{deg}(y_i)=(1,-1)$ proves (c). It remains to show that (a) implies (b) and that (d) implies (a). 

\subsection{} For each partition $\lambda$ of $n$, let $S^\lambda$ be the corresponding irreducible $S_n$-module.  Define the \emph{standard module} $$\Delta(\lambda)=\mathrm{Ind}^H_{\CC[\hh^*] \rtimes S_n} S^\lambda$$ and the \emph{baby standard module} (termed ``baby Verma" module by Iain Gordon \cite{Gor}, in analogy with Lie theory)
$$\overline{\Delta}(\lambda)=\Delta(\lambda) / I \Delta(\lambda),$$ where $I$ is the ideal in $\CC[\hh]$ generated by the positive degree $S_n$ invariants. 

\subsection{} Inspection of the defining relations shows that there are an involutive automorphism $\phi$ and an involutive anti-automorphism $\iota$ of $H$ given by $\phi(x_i) = y_i=\iota(x_i)$, $\phi(y_i)=x_i=\iota(y_i)$, $\phi(w)=\mathrm{det}(w) w$, and $\iota(w)=w^{-1}$ for $w \in S_n$. These commute, and their composition is an involutive anti-automorphism fixing $x_i$'s and $y_i$'s and mapping $w$ to $\mathrm{det}(w) w^{-1}$.  Given an $H$ module $M$, we define the twist $M^\phi$ of $M$ by $\phi$ to be the same as a vector space but with a new action $\cdot$ given by the formula $f \cdot m=\phi(f) m$. Then $L(\lambda)^\phi$ is irreducible, hence isomorphic to some $L(\mu)$.  The automorphism $\phi$ is crucial for our proof that (d) implies (a).

\subsection{} There is a \emph{contravariant form} $\la \cdot,\cdot \ra'$ on $\Delta(\lambda)$ such that 
$$L(\lambda)=\Delta(\lambda) / \mathrm{Rad} (\la \cdot,\cdot \ra ').$$ It is defined by the requirements that it be $\CC$-bilinear, restrict to a fixed non-degenerate $W$-invariant bilinear form on the degree zero part $S^\lambda \subseteq L(\lambda)$, and 
\begin{equation} \label{contravar}
\la f g,h \ra'=\la g, \iota(f) h \ra'
\end{equation} for all $f \in H$ and $g,h \in \Delta(\lambda)$. Evidently the contravariant form descends to a non-degenerate pairing on $L(\lambda)$. Furthermore, since distinct graded pieces of $\Delta(\lambda)$ are orthogonal for the contravariant form, its radical is a graded subspace and hence $L(\lambda)$ is a graded quotient. We write $L(\lambda)^d$ for the degree $d$ piece.

\subsection{} \label{single variable} It follows from \cite{Gor} or the norm calculation carried out by Dunkl in \cite{Dun} that the image of the minimal occurence of the trivial representation in $\CC[x_1,\dots,x_n] \otimes S^\lambda$ has non-zero image $t_\lambda$ in $L(\lambda)$. By irreducibility of $L(\lambda)$, 
$$\{f \in L(\lambda) \ | \ yf=0 \ \hbox{for all $y \in \hh$} \}=L(\lambda)^{-n(\lambda)}=S^\lambda$$ is the degree $-n(\lambda)$ piece of $L(\lambda)$. Therefore the $\CC[\hh^*]$-submodle $\CC[\hh^*] t_\lambda \subseteq L(\lambda)$ contains $S^\lambda$ as its socle, and hence the kernel $I$ of the map $\CC[\hh^*] \rightarrow \CC[\hh^*] t_\lambda$ does not contain the copy of $S^\lambda$ in degree $n(\lambda)$. It follows that $I \subseteq I_\lambda$. But since $S^\lambda$ is the socle of $\CC[\hh^*] t_\lambda$ and the natural map $\CC[\hh^*] t_\lambda \cong \CC[\hh^*]/I \rightarrow \CC[\hh^*]/I_\lambda=R_\lambda(y)$ does not kill $S^\lambda$, it must be an isomorphism (this fact was also proved by Gordon in the proof of Theorem 6.7, \cite{Gor}). 

\subsection{} Define $$z_i=y_i x_i+\sum_{1 \leq j <i} s_{ij}.$$ By Theorem 5.1 of \cite{Gri} each joint eigenvalue of $z_1,\dots,z_n$ on $L(\lambda)$ is a permutation of the numbers $-\mathrm{ct}(b)$ for boxes $b \in \lambda$. Since $$\phi(z_i)=x_i y_i-\sum_{1 \leq j < i} s_{ij}=y_i x_i+\sum_{i<j \leq n} s_{ij}=w_0 z_{w_0(i)} w_0^{-1},$$ the eigenvalues of $z_i$ on $L(\lambda)^\phi$ are the same. Thus since the set $\{\mathrm{ct}(b) \ | \ b \in \lambda \}$ determines $\lambda$ it follows that $L(\lambda)^\phi \cong L(\lambda)$. This gives a linear isomorphism $\psi: L(\lambda) \rightarrow L(\lambda)$ such that $\psi(h f)=\phi(h) \psi(f)$ for $h \in H$ and $f \in L(\lambda)$. Since $\phi^2=1$, the map $\psi^2$ is an $H$-module endomorphism, hence equal to a scalar. By renormalizing we will assume that this scalar is $1$. 

\subsection{} Note that $d_\lambda=\psi(t_\lambda)$ is an  occurence of the determinant representation in $L(\lambda)$. Again, the norm calculation of \cite{Dun} or the paper \cite{Gor} shows that $d_\lambda$ is in degree $n(\lambda^t)$ of $L(\lambda)$, and now the same argument as in \ref{single variable} proves that the space $\CC[\hh^*] d_\lambda$ is isomorphic to $R_{\lambda^t}(y) \otimes \mathrm{det}$. Applying $\psi$ to this shows that the space $\CC[\hh] t_\lambda$ is isomorphic to $R_{\lambda^t}(x)$. By irreducibility we have $L(\lambda)=\CC[\hh] \CC[\hh^*] t_\lambda=\CC[\hh^*] \CC[\hh] t_\lambda$, and thus $\mathrm{gr} L(\lambda)$ is a quotient of $R_{\lambda^t}(x) \otimes_\CC R_\lambda(y)$ in which the sign representation survives.

\subsection{} By Proposition 4.1.2 of \cite{Hai2}, the kernel $J_\lambda$ of the quotient map $\CC[x_1,\dots,x_n,y_1,\dots,y_n] \rightarrow R_\lambda$ consists of those $f$ such that the principal ideal generated by the image $\overline{f}$ of $f$ in $R_{\lambda^t}(x) \otimes_\CC R_\lambda(y)$ has zero intersection with the unique copy of the sign representation. If $I_\lambda$ is the kernel of the quotient map $\CC[x_1,\dots,x_n,y_1,\dots,y_n] \rightarrow \mathrm{gr}(L(\lambda))$, then every $f \in I_\lambda$ must have this property, since the sign representation survives in $\mathrm{gr} L(\lambda)$. Thus  $I_\lambda \subseteq J_\lambda$ and $R_\lambda$ is a quotient of $\mathrm{gr} L(\lambda)$. Therefore (a) implies (b).

\subsection{}  Next we prove that (d) implies (a). From now on we will abbreviate by omitting $\lambda$ when appropriate, in particular writing $L=L(\lambda)$. We note that it suffices to prove that $\mathrm{gr} L$ is Gorenstein, for then the quotient map $\mathrm{gr} L \rightarrow R_\lambda$ constructed in the previous section will necessarily be an isomorphism. To prove that $\mathrm{gr} L$ is Gorenstein, it suffices to construct a non-degenerate bilinear form on it, compatible with the operators $x_i$ and $y_i$.

\subsection{} Define a non-degenerate pairing 
$$\la \cdot,\cdot \ra : L \times L \rightarrow \CC$$ by the formula 
$$\la f,g \ra=\la \psi(f),g \ra' \quad \hbox{for $f,g \in L$.}$$ This is compatible with the $H$-action in the sense that 
\begin{equation} \label{compat} 
\la h f_1,f_2 \ra=\la f_1, \iota \phi(h) f_2 \ra
\end{equation} for $h \in H$ and $f_1,f_2 \in L$. In particular,
\begin{equation} \label{compat2}
\la x_i f_1,f_2 \ra=\la f_1,x_i f_2 \ra \quad \text{and} \quad \la y_i f_1,f_2 \ra=\la f_1,y_i f_2 \ra \quad \hbox{for all $1 \leq i \leq n$ and $f_1,f_2 \in L$.}
\end{equation}

\subsection{} Here we prove that $\la \cdot,\cdot \ra$ is either symmetric or skew-symmetric. Observe first that $\la \psi(\cdot),\psi(\cdot) \ra'$ is a non-degenerate bilinear form on $L$ that also satisfies the contravariance \eqref{contravar}.Therefore there is some $c \in \CC$ with $\la \psi(\cdot),\psi(\cdot) \ra'=c \la \cdot,\cdot \ra'$. Now compute
$$\la f,g \ra=\la \psi(f),g \ra'=\la g,\psi(f) \ra'=c \la \psi(g),f \ra'=c \la g,f \ra.$$ Repeating shows $c^2=1$, as was to be proved. From now on, for a subspace $VÊ\subseteq L$ we will employ the notation
$$V^\perp=\{f \in L \ | \ \la f, V \ra=0 \}=\{f \in L \ | \ \la V, f \ra=0 \}.$$

\subsection{} Now we show that the form $\la \cdot,\cdot \ra$ descends to $\mathrm{gr} L$. First, observe that $S_n$-equivariance implies $\la t_\lambda, f \ra=0$ for all $f \in L$ of isotype distinct from $\mathrm{det}$. Since the determinant occurs only in bidegree $(n(\lambda^t),n(\lambda))$ in $R_\lambda$, it follows from the surjection that we have constructed $\mathrm{gr} L \rightarrow R_\lambda$ that the unique occurence of $\mathrm{det}$ in $L$ does not occur in $W_l L$ for $l <  n(\lambda^t)+n(\lambda)$. Hence $\la t_\lambda, W_l L \ra=0$ unless $l \geq n(\lambda^t)+n(\lambda)$. Using the relation \eqref{compat} and induction we obtain
\begin{equation} \label{orth} 
\la W_l L , W_{l'} L \ra = 0 \quad \hbox{unless $l+ l' \geq n(\lambda^t)+n(\lambda)$.}
\end{equation} Define a pairing $\mathrm{gr} \la \cdot,\cdot \ra$ on $\mathrm{gr} L $ by declaring $\mathrm{gr}^l L$ orthogonal to $\mathrm{gr}^{l'} L$ unless $l+l'=n(\lambda^t)+n(\lambda)$, in which case we define
$$\mathrm{gr} \la \overline{f_1},\overline{f_2} \ra=\la f_1,f_2 \ra \quad \hbox{for $f_1 \in W_l L$ and $f_2 \in W_{l'}L$,}$$ which is well-defined thanks to \eqref{orth}. By definition it satisfies the analog of \eqref{compat2}.

\subsection{} Write $m=n(\lambda^t)+n(\lambda)$ and, for $0 \leq l \leq m$, write $l'=m-l$. Noting that $W_{l'-1} L \subseteq W_l L ^ \perp$, the non-degeneracy of $\mathrm{gr} \la \cdot,\cdot \ra$ is equivalent to the equalities
\begin{equation*}
\mathrm{dim}(W_{l'-1} L)=\mathrm{codim}(W_{l}L) \quad \hbox{for $0 \leq l \leq m$.}
\end{equation*} In turn these are equivalent to
\begin{equation} \label{dims}
\mathrm{dim}(\mathrm{gr}^{l'} L)=\mathrm{dim}(\mathrm{gr}^l L) \quad \hbox{for $0 \leq l \leq m$.}
\end{equation} This completes the proof that (d) implies (a) and hence the proof of the theorem.

\subsection{} If we could establish the analogue of the final isomorphism of Deligne's paper \cite{Del}, we would have a proof of the equalities \eqref{dims} and hence the $n!$ theorem. Perhaps the techniques of \cite{VaVa} could be adapted to identify the module $L$ with the cohomology ring of an affine Springer fiber.

\subsection{Wreath product groups} Using the norm calculation from \cite{DuGr}, one could define a version of Macdonald polynomials for wreath product groups $G(r,1,n)$ in an analogous fashion. At the moment we do not see a way to connect these to the polynomials defined in \cite{Hai2} and studied by Bezrukavnikov-Finkelberg in \cite{BeFi}, or those studied in \cite{BDLM1} and \cite{BDLM2}. In any case, one should hope for a reasonable formula connecting Macdonald polynomials and the bigraded character of the quotient of the coinvariant ring studied in \cite{GoGr}.

\def\cprime{$'$} \def\cprime{$'$}


\begin{thebibliography}{GGOR}

\bibitem[As]{As} S. Assaf, \emph{Dual equivalence graphs I: A combinatorial proof of LLT and Macdonald positivity}, 
arXiv:1005.3759

\bibitem[BDLM1]{BDLM1} O. Blondeau-Fournier, P. Desrosiers, L. Lapointe, P. Mathieu, \emph{Macdonald polynomials in superspace: conjectural definition and positivity conjectures},  arXiv:1112.5188

\bibitem[BDLM2]{BDLM2} O. Blondeau-Fournier, P. Desrosiers, L. Lapointe, P. Mathieu, \emph{Macdonald polynomials in superspace as eigenfunctions of commuting operators},  arXiv:1202.3922

\bibitem[BeFi]{BeFi} R. Bezrukavnikov and M. Finkelberg, \emph{Wreath Macdonald polynomials and categorical McKay correspondence} arXiv:1208.3696

\bibitem[Del]{Del} P. Deligne, \emph{Th\'eorie de Hodge I}, Actes, Congr\`es intern. math., 1970. Tome 1, p. 425--430

\bibitem[Dun]{Dun} C. Dunkl, \emph{Symmetric and antisymmetric vector-valued Jack polynomials}, Sem. Lothar. Combin. B64a (2010), 31 p., arXiv:1001.4485

\bibitem[DuGr]{DuGr} C. Dunkl and S. Griffeth, \emph{Generalized Jack polynomials and the representation theory of rational Cherednik algebras}, Selecta Math. (N.S.) 16 (2010), 791-818,  arXiv:1002.4607

\bibitem[EtGi]{EtGi} P. Etingof and V. Ginzburg, \emph{Symplectic reßection algebras, Calogero-Moser space, 
and deformed Harish-Chandra homomorphism}, Invent. Math. 147 (2002), no. 2, 243-348, arXiv:math/0011114

\bibitem[GaHa]{GaHa} A. Garsia and M. Haiman, \emph{A graded representation model for MacdonaldÕs polynomials}, Proc. Nat. Acad. 
Sci. U.S.A. 90 (1993), no. 8, 3607Ð3610. 

\bibitem[GaPr]{GaPr} A. Garsia and C. Procesi, \emph{On certain graded $S_n$ -modules and the $q$-Kostka polynomials}, 
Adv. Math. 94 (1992), no. 1, 82Ð138

\bibitem[Gor]{Gor} I. Gordon, \emph{Baby Verma modules for rational Cherednik algebras}, Bulletin of the London Mathematical Society (2003), 35, pp 321-336, arXiv:math/0202301

\bibitem[Gor2]{Gor2} I. Gordon, \emph{On the quotient ring by diagonal invariants}, Invent. Math.153, no. 3, (2003), 503-518

\bibitem[Gor3]{Gor3} I. Gordon, \emph{Macdonald positivity via the Harish-Chandra $D$-module}, to appear Invent. Math.

\bibitem[GoGr]{GoGr} I. Gordon and S. Griffeth, \emph{Catalan numbers for complex reflection groups}, to appear Amer. Jour. Math.

\bibitem[Gri]{Gri} S. Griffeth, \emph{Orthogonal functions generalizing Jack polynomials}, Trans. Amer. Math. Soc. 362 (2010), 6131-6157

\bibitem[GrHa]{GrHa} I. Grojnowski and M. Haiman, \emph{Affine Hecke algebras and positivity of LLT and Macdonald polynomials}, 
preprint (2007)

\bibitem[Hai]{Hai} M. Haiman, \emph{Hilbert schemes, polygraphs and the Macdonald positivity conjecture}, J. Amer. Math. Soc. 14 (2001), 
no. 4, 941Ð1006, arXiv:math.AG/0010246

\bibitem[Hai2]{Hai2} M. Haiman, \emph{Combinatorics, symmetric functions, and Hilbert schemes}, Current Developments in Mathematics 2002, no. 1 (2002), 39-111

\bibitem[Hai3]{Hai3} M. Haiman, \emph{Macdonald polynomials and geometry}, New perspectives in geometric combinatorics 
(Billera, Bj¬orner, Greene, Simion, and Stanley, eds.), MSRI Publications, vol. 38, Cambridge 
University Press, 1999, pp. 207Ð254

\bibitem[HoNa]{HoNa} R. Holmes and D. Nakano, \emph{Brauer-type reciprocity for a class of graded associative algebras}, J. Algebra 144 (1991), 117--126

\bibitem[KaRo]{KaRo} M. Kashiwara and R. Rouquier, \emph{Microlocalization of rational Cherednik algebras}, Duke Math. J. 144 (2008), no. 3, 525Ð573, arXiv:0705.1245

\bibitem[Mac]{Mac} I.G. Macdonald, \emph{Symmetric functions and Hall polynomials}, Second edition. With contributions by A. Zelevinsky. Oxford Mathematical Monographs. Oxford Science Publications. The Clarendon Press, Oxford University Press, New York, 1995

\bibitem[VaVa]{VaVa} M. Varagnolo and E. Vasserot, \emph{Finite-dimensional representations of DAHA and affine Springer fibers: the spherical case}, Duke Math. J. 147 (2009), no. 3, 439Ð540, arXiv:0705.2691
\end{thebibliography}
\end{document}